\documentclass[titlepage,14pt]{extarticle}
\usepackage{amsfonts,amsmath,amssymb,indentfirst}

\usepackage{geometry}
\def\inbr#1{ \{ #1 \} }

\def\N{{\mathbb N}}
\def\Q{{\mathbb Q}}
\def\B{{\mathcal{B}}}
\def\M{{\mathcal{M}}}
\def\D{{\mathbb D}}

\newtheorem{lemma}{Lemma}
\newtheorem{corollary}{Corollary}
\newtheorem{conjecture}{Conjecture}
\def\proof{{\bf Proof. }}

\allowdisplaybreaks[4]

\begin{document}

\begin{center}
{\bf \scshape A Note on Arithmetical Properties of Multiple Zeta Values} \\
{\bf \scshape Sergey Zlobin}
\end{center}

\vskip 0.5cm

Multiple zeta values
$$
\zeta(s_1, s_2, \dots, s_l) = \sum_{n_1 > n_2 > \cdots > n_l \ge 1}
\frac{1}{n_1^{s_1} \cdots n_l^{s_l}}, \quad s_j \in \N
$$
are actively studied,
but the majority of results are various
identities between these values.
In this paper we touch their arithmetical properties.

The {\it weight} of the vector $\vec{s}=(s_1, s_2, \dots, s_l)$
(or of the series $\zeta(s_1, s_2, \dots, s_l)$)
is the sum $s_1+s_2+\dots+s_l$. Let $\D_w$ be the dimension (over $\Q$) of the 
$\Q$-vector space generated by all multiple zeta values of weight
$w \ge 2$.

The following conjecture on $\D_w$ appears in \cite{zagier}. 
\begin{conjecture}
\label{con1}
If $w \ge 2$, then $\D_w=d_w$, where the sequence $\{ d_w \}$
is given by the recurrence
$d_w=d_{w-3}+d_{w-2}$ with initial values $d_0=1$, $d_1=0$, $d_2=1$.
\end{conjecture}
The numbers $d_w$ can be defined by the generating function
$$
\sum_{w=0}^{\infty} d_w x^w = \frac{1}{1-x^2-x^3}.
$$
By Conjecture \ref{con1} the number of linearly independent numbers among
the $2^{w-2}$ multiple zeta values of weight $w$ is
$d_w = O(\alpha^w)$ as $w \to \infty$, where $\alpha=1.32..$ is the 
root of $x^3-x-1=0$.
So there should be a lot of
(independent) relations among these values. Having certain relations,
we can prove an upper estimate for $\D_w$. The paper \cite{terasoma}
claims to prove the inequality $\D_w \le d_w$. However,
the proof is very complicated.

A lower estimate for $\D_w$ is connected with the arithmetical nature of 
multiple zeta values. With the help of multiplication by $\zeta(w_1-w_2)$,
it is easy to show 
that $\D_{w_1} \ge \D_{w_2}$ if $w_1 \ge w_2+2$. However, it is not proved
that $\D_w > 1$ for a single $w$!

Choose the
following sets among all vectors with positive integer components:
$$
\B=\{ \vec{s} : s_i \in \{ 2, 3 \} \}, \quad
\B_w = \{ \vec{s} \in \B : w( \vec{s} ) = w \}.
$$
M.~Hoffman (\cite{hoffman}) made the following conjectures.

\begin{conjecture}
\label{con2}
For any $\vec{s_0}$ the value
$\zeta(\vec{s_0})$ can be represented 
as a linear form with
rational coefficients in values $\zeta(\vec{s})$ with $\vec{s} \in
\B_{w(\vec{s_0})}$.
\end{conjecture}

This conjecture has been checked by Hoang Ngoc Minh
for $\vec{s_0}$ with weight $\le 16$.

\begin{conjecture}
\label{con3}
If $\vec{s_1}$, \dots, $\vec{s_n}$ are distinct elements of $\B$, then
1, $\zeta(\vec{s_1})$, \dots, $\zeta(\vec{s_n})$ are
linearly independent over $\Q$.
\end{conjecture}

Suppose that Conjecture \ref{con3} is true; then the linear form representation
from Conjecture \ref{con2} is unique.
From Conjectures \ref{con2} and \ref{con3}
there follows Conjecture \ref{con1} (see \cite{hoffman}).

Since $\zeta{(\inbr{2}_k)} = \frac{\pi^{2k}}{(2k+1)!}$
(see \cite[(36)]{bradley}), these
values are irrational (moreover they and 1 are linearly
independent over $\Q$). Also, by Apery's theorem, the number $\zeta(3)$
is irrational.
There is no certainty about the arithmetic nature of $\zeta(\vec{s})$ for 
any other $\vec{s} \in \B$.

Let $\zeta(\vec{s_0}) \in \Q$ and $w(\vec{s_0})$ be odd.
Suppose that $\zeta(\vec{s_0}) \zeta(2k)$ is represented as 
a linear combination with rational coefficients in numbers $\zeta(\vec{s})$,
$\vec{s} \in \B_{w(\vec{s_0})+2k}$ (as implied by Conjecture \ref{con2});
then  there is at least one irrational among these numbers.
For instance if $\zeta(2,3) \in \Q$ or
$\zeta(3,2) \in \Q$ then one of the numbers $\zeta(3,2,2)$,
$\zeta(2,3,2)$ and $\zeta(2,2,3)$ is irrational. Similarly,
suppose that $\zeta(\vec{s_0}) \in \Q$, $w(\vec{s_0})$ is even and
$\zeta(\vec{s_0}) \zeta(3)$ is represented as a linear
combination with rational coefficients in numbers $\zeta(\vec{s})$,
$\vec{s} \in \B_{w(\vec{s_0})+3}$; then there is at
least one irrational among them.

Furthermore, we prove a certain result about the linear
independence of multiple zeta values of different weights.

\begin{lemma}
\label{GenLinIndepLemma}
Let $x \notin \Q$, numbers $y_i$, $i=1,
\dots, k$ be such that 1, $y_1$, \dots, $y_k$ are linearly independent
over $\Q$. Then there exist $k-1$ numbers among the $x y_i$ such that
1, $x$ and they are linearly independent over $\Q$.
\end{lemma}
\proof
We prove by contradiction. Let the numbers 1, $x$, $x y_1$, \dots, $x y_{k-1}$
be linearly dependent over $\Q$. I.e. there exist
integers $A_1$, $B_1$ and $C_{1i}$ not all zero
such that
$$
A_1+B_1 x+ \sum_{i=1}^{k-1} C_{1i} x y_i =0.
$$
If $A_1=0$ then divide this equality by $x$; we get that
1, $y_1$, \dots, $y_{k-1}$ are linearly dependent, which contradicts
the hypothesis. If all $C_{1i}=0$ then $x$ is rational.
Hence there exists $p \in [1,k-1]$ such that $C_{1p} \ne 0$. Let
integers $A_2$, $B_2$ and $C_{2i}$ be not all zero
such that
$$
A_2+B_2 x+ \sum_{1 \le i \le k, i \ne p} C_{2i} x y_i =0.
$$
Similarly $A_2 \ne 0$. Multiply the first equality by
$A_2$ and subtract the second equality, multiplied by $A_1$.
We get (letting $C_{1k}=0$, $C_{2p}=0$)
$$
(B_1 A_2 - B_2 A_1)x + \sum_{i=1}^{k} (C_{1i} A_2 - C_{2i} A_1) x y_i =0.
$$
Divide this equality by $x$. Then we get a linear form
in 1, $y_1$, \dots, $y_k$;
moreover the coefficient of $y_p$ is $C_{1p} A_2 \ne 0$,
which contradicts the linear independence of 1 and the numbers $y_i$.
The lemma is proved.

\begin{corollary}
For any  positive integer $l$, the numbers 1, $\zeta(3)$
and certain $l$ numbers
from $\zeta(3) \zeta(2k)$, $k=1,\dots,l+1$ are linearly independent over
$\Q$.
\end{corollary}
\proof
Take $x=\zeta(3)$, $y_k=\zeta(2k)$ in Lemma $\ref{GenLinIndepLemma}$.

This corollary and the equality
$$
\zeta(3) \zeta(2k) = \zeta(2k+3)+ \zeta(3,2k) + \zeta(2k,3)
$$
yield another
\begin{corollary}
\label{MultyZetaLinIndepCor}
For each $w \ge 2$, let $\M_w$ be a set of vectors of
weight $w$ such that all multiple zeta values of the weight $w$
are rational linear combinations of $\zeta(\vec{s})$ with
$\vec{s} \in \M_w$. Then for any positive integer $l$,
there exists a subset $I$ of $\{5,7,\dots,2l+5 \}$ with $\#I=l$, and
$l$ vectors $\vec{t_i} \in \M_i$,  $i$ ranging over $I$,
such that 1, $\zeta(3)$ and the numbers
$\zeta(\vec{t_i})$ are linearly independent over $\Q$.
\end{corollary}
From the equality $\zeta{(\inbr{2}_k)} = \frac{\pi^{2k}}{(2k+1)!}$
it is clear that
$$
\dim_{\Q} (\Q \oplus \bigoplus_{\vec{s} \in
\B_2 \cup \B_4 \cup \cdots \cup \B_{2l}}
\Q \zeta(\vec{s}) ) \ge l+1.
$$
By Conjecture \ref{con2} it is possible to take
$\B_w$ as the set $\M_w$ required in Corollary \ref{MultyZetaLinIndepCor}.
If so, then
$$
\dim_{\Q} (\Q \oplus \bigoplus_{\vec{s} \in
\B_3 \cup \B_5 \cup \cdots \cup \B_{2l+5}}
\Q \zeta(\vec{s}) ) \ge l+2.
$$
Since Conjecture \ref{con2} has been checked for vectors of weights $\le 16$,
this estimate is not conditional for $l \le 5$.

\begin{corollary}
\label{TwoThreeCor}
There exists
$$\vec{s_0} \in \{(2,3), (3,2), (2,2,3), (2,3,2), (3,2,2) \},$$
such that the numbers 1, $\zeta(3)$ and $\zeta(\vec{s_0})$
are linearly independent over $\Q$.
\end{corollary}
\proof We use Corollary \ref{MultyZetaLinIndepCor} with
$l=1$ choosing $M_5=\{ (2,3), (3,2) \}$ and $M_7=\{ (2,2,3),
(2,3,2), (3,2,2) \}$.

We emphasize that the result of Corollary \ref{TwoThreeCor} is unconditional.

\vskip 3mm

I thank kindly J. Sondow for several suggestions.


\newcommand{\namefont}{\scshape}
\newcommand{\titlefont}{\itshape}
\def\nomer{No.}

\end {document}